\documentclass[10pt,letterpaper]{article}%
\usepackage[utf8]{inputenc}
\usepackage{amsmath}
\usepackage{amsfonts}
\usepackage{amssymb}
\usepackage{graphicx}%
\providecommand{\U}[1]{\protect\rule{.1in}{.1in}}
\newtheorem{theorem}{Theorem}

\newtheorem{corollary}[theorem]{Corollary}

\newtheorem{lemma}[theorem]{Lemma}

\newenvironment{proof}[1][Proof]{\noindent\textbf{#1.} }{\ \rule{0.5em}{0.5em}}
\begin{document}

\title{Existence of periodic solutions for a scalar differential equation modelling
optical conveyor belts}
\author{Luis Carretero$^{1}$, Jos\'{e} Valero$^{2}$\\$^{1}${\small Departamento de Ciencia de Materiales, \'{O}ptica y
Tecnolog\'{\i}a Electr\'{o}nica,}\\{\small Universidad Miguel Hern\'{a}ndez de Elche, 03202, Elche, Spain}\\$^{2}${\small Centro de Investigaci\'{o}n Operativa,}\\{\small Universidad Miguel Hern\'{a}ndez de Elche, 03202, Elche, Spain}}
\date{}
\maketitle

\begin{abstract}
We study a one-dimensional ordinary differential equation modelling optical
conveyor belts, showing in particular cases of physical interest that periodic
solutions exist. Moreover, under rather general assumptions it is proved that
the set of periodic solutions is bounded.

\end{abstract}

\bigskip


\textbf{Subject Mathematics Classification (2010): }\ 78A10, 34C25, 34B15

\section{Introduction}

Electromagnetic fields can exert forces on micro and nano-particles which are
the result of radiation scattering produced by the own particle. As it was
demonstrated by Ashkin \cite{Ashkin1970}, these electromagnetic forces can be
used to trap and manipulate the matter. Since then, electromagnetic
microparticle manipulation has become an important technique in a wide range
of fields like biology, colloid science or microfluidics.

In the case of small particles (radius of a few nanometers), assuming that we
are working far from the resonance of nano-particles region, the main optical
techniques for micromanipulation like optical tweezers, tractor beams or
conveyor belts are based on the gradient forces
\cite{Chaumet00,Vesperinas2010,Ruffner2012,Carretero2014} generated by the
interaction of the particles with an spatially inhomogeneous optical beams
like Gaussian beams or Bessel beams \cite{Saleh91}.

In this work, we are going to study the particular case of an optical conveyor
belt \cite{Cizmar2005,Cizmar2006}, where the gradient forces that act on the
z-axially confined particles are obtained by means of the superposition of two
temporally dephased counter propagating complex electromagnetic fields of
frequency $\omega$. Thus, by assuming linearly polarized beams along the
$x$-axis, the electric field is given by:
\begin{equation}
\vec{E}(t,z)=\big(e_{0}f_{0}(z)\exp(ikz)+e_{0}f_{0}(z)\exp(-ikz)\exp
(ibt)\big)\exp(i\omega t), \label{campo}%
\end{equation}
where $k=n\omega/c$ is the wavenumber in a medium with refractive index $n$,
$b$ is the parameter that makes the conveyor work by controlling the relative
phase of counter propagating fields \cite{Ruffner2012}, $e_{0}$ is the
amplitude of the electric field and $f_{0}(z)$ takes into account the spatial
variation of the electric field amplitude.

If particles are small enough (Rayleigh regimen), the matter-radiation
interaction can be modelled by electric dipoles \cite{Chaumet00}:
\begin{equation}
\vec{p}=4\pi\epsilon_{0}\alpha\vec{E},
\end{equation}
where $\epsilon_{0}$ is the permittivity of the free space and $\alpha$ is the
polarizability of the particle.

The dipole particle dynamic inside the media is governed by the differential
equation:
\begin{equation}
m\,z^{\prime\prime}(t)=-\gamma\,z^{\prime}(t)+F(t,z) \label{eq1}%
\end{equation}
where $\gamma$ is the friction coefficient of the medium where, particles of
mass $m$ are immersed, and $F(t,z)$ correspond to the time average axial force
that acts on a Rayleigh particle inside of the electromagnetic field
\eqref{campo}, that it is given by \cite{Chaumet00}:
\begin{equation}
F(t,z)=4\pi\epsilon_{0}\mathcal{R}[\alpha]\frac{\partial|\vec{E}(t,z)|^{2}%
}{\partial z}, \label{fuerza}%
\end{equation}
where $\mathcal{R}[$\textperiodcentered$]$ denote the real part. For obtaining
the last expression, as it has been previously mentioned, we have neglected
the imaginary part of dipole moment assuming that the working wavelength of
the electromagnetic fields are far from the particles resonance.

For Rayleigh particles, differential equation \eqref{eq1} is over-damped
\cite{Marago2013} ($m\,z^{\prime\prime}(t)\approx0$), so the particle dynamic
inside the media of refractive index $n$ is governed by the non autonomous
differential equation:
\begin{equation}
z^{\prime}=\frac{F(t,z)}{\gamma}=F_{z}(t,z). \label{langevin}%
\end{equation}
Introducing equation \eqref{campo} in \eqref{fuerza} and \eqref{langevin}, it
can be observed that the force $F_{z}$ is a gradient one, i.e., $F_{z}%
=\frac{\partial V(t,z)}{\partial z}$, and then the particles in this field are
immersed in a potential energy given by:
\begin{equation}
V(t,z)=\frac{4\pi\epsilon_{0}\mathcal{R}[\alpha]}{\gamma}|\vec{E}|^{2}%
=F_{0}\,f(z)\cos\left(  kz-\frac{bt}{2}\right)  ^{2} \label{potencial}%
\end{equation}
where we have introduced the constant $F_{0}=\frac{4\pi\epsilon_{0}%
\mathcal{R}[\alpha]e_{0}^{2}}{\gamma}$ and the function $f(z)=|f_{0}(z)|^{2}$.
Constant $F_{0}$ is proportional to total intensity of the electromagnetic
wave. The periodic function that depends on the phase difference of the
electromagnetic field give rise to the optical conveyor, whose continuous
variation move the trapped particles along the $z$-axis.

Then, the dynamic of particles in an axial optical conveyor belt can be
modeled by the differential equation:
\begin{equation}
z^{\prime}=F_{z}(t,z), \label{problem}%
\end{equation}
where%
\begin{align}
F_{z}(t,z)  &  =\frac{\partial V(t,z)}{\partial z},\label{Funcion}\\
V(t,z)  &  =F_{0}\,f(z)\cos\left(  kz-\frac{bt}{2}\right)  ^{2},
\label{Potencial}%
\end{align}
with $F_{0},k,b>0$.

As we have mentioned, the potential $V(t,z)$ is directly proportional to the
conveyor's axial intensity, so the $f(z)$ function determines the axial region
where the conveyor belt has enough strength to move the particles. We are
going to analyze three different cases for the function $f(z)$:

\begin{enumerate}
\item Constant axial region strength $f(z)=1.$

In this case the conveyor strength is not spatially limited and can be
physically obtained by means of the interference of two dephased counter
propagating plane waves.

\item Lorentzian axial region strength $f(z)=\frac{1}{1+(z/z_{0})^{2}}%
,\ z_{0}>0.$

This behavior can be physically obtained through the interference of two
dephased counter propagating Gaussian beams. In this case the parameter
$z_{0}$ is known as the Rayleigh range of the Gaussian beams, and the axial
region strength has its maximum at $z=0$ position and drops gradually as
$\left\vert z\right\vert $ increases, reaching half its peak value at $z=\pm
z_{0}$, i.e. $(f(\pm z_{0})=\frac{1}{2})$ \cite{Saleh91}.

\item Gaussian axial region strength $f(z)=\exp(-2\frac{z^{2}}{z_{0}^{2}%
}),\ z_{0}>0.$

In this case we have also the axial region strength at $z=0$ and $f$ drops
monotonically as $\left\vert z\right\vert $ increases. The meaning of $z_{0}$
is different to the previous case since $f(z_{0})=\exp(-2)$. Also, it is more
difficult to obtain physically this kind of conveyor belt than previous ones;
however, similar methodologies to those described in \cite{Zamboni2004,
Zamboni2005} could be used in order to obtain it.
\end{enumerate}

The main goal of this paper is to establish the existence of periodic
solutions for equation \eqref{problem} in the cases of Lorentzian and Gaussian
axial region strengths. This result was suggested in \cite{Carretero2015} by
numerical simulations for a four-dimensional model, and now we give a rigorous
mathematical proof of this statement in a simpler one-dimensional model of
optical conveyor belts.

This type of problem is interesting for both physical and mathematical points
of view, as the existence of periodic solutions for equations of this kind is
not known in the mathematical literature as far as we know. In order to prove
it we have used an abstract result about the existence of solutions of
boundary value problems given in \cite{Kiguradze97}. We observe that standard
techniques such as obtaining an invariant region in order to use the fixed
point theorem (see e.g. \cite{Cronin, ShiSongLi, Zhao}), the use of a Lyapunov
function \cite{Xiang} or degree theory \cite{Amann} seem to be uselees for our
problem. Thus, a very specific argument has been developed for our particular
type of equations.

\section{Non-plane waves: $f(z)\neq1$}

Let us consider the differential equation%
\begin{equation}
z^{\prime}=F_{z}(t,z) \label{Equation}%
\end{equation}
with $F_{z}$ defined by \eqref{Funcion}--\eqref{Potencial} and with $f\left(
z\right)  $ not identically equal to $1$. In fact, $f\left(  z\right)  $ will
be a function satisfying the asymptotic behaviour $f\left(  z\right)
\rightarrow0$ as $\left\vert z\right\vert \rightarrow\infty.$

First, we prove that under rather general assumptions on the function $f$ the
set of periodic solutions with period $T=\frac{4\pi}{b}$ is bounded.

Second, we obtain the existence of periodic solutions with period
$T=\frac{4\pi}{b}$ for the particular cases $f(z)=\frac{1}{1+(z/z_{0})^{2}}$
and $f(z)=exp(-2\frac{z^{2}}{z_{0}^{2}})$.

\subsection{Boundedness of the set of periodic solutions}

We assume that $f$ satisfies the following conditions:

\begin{enumerate}
\item $f\in C^{1}\left(  \mathbb{R}\right)  .$

\item $f\left(  z\right)  \rightarrow0$ if $\left\vert z\right\vert
\rightarrow\infty.$

\item $f^{\prime}\left(  z\right)  \rightarrow0$ if $\left\vert z\right\vert
\rightarrow\infty.$

\item If $z_{n},v_{n}\rightarrow\infty$, $\frac{v_{n}}{z_{n}}\rightarrow1$, as
$n\rightarrow\infty$, and $f^{\prime}\left(  z_{n}\right)  \not =0$, then%
\begin{equation}
\frac{f^{2}\left(  v_{n}\right)  }{f^{\prime}\left(  z_{n}\right)
}\rightarrow0,\ \frac{\left(  f^{\prime}\left(  v_{n}\right)  \right)  ^{2}%
}{f^{\prime}\left(  z_{n}\right)  }\rightarrow0. \label{Cond4Functf}%
\end{equation}

\item The flow generated by equation \eqref{Equation} has no fixed points.
\end{enumerate}

We note that%
\[
F_{z}(t,z)=-kF_{0}f\left(  z\right)  \sin\left(  2kz-bt\right)  +F_{0}\left(
\cos\left(  kz-\frac{bt}{2}\right)  \right)  ^{2}f^{\prime}\left(  z\right)
.
\]

The function $t\mapsto F_{z}\left(  t,z\right)  $ is periodic with period
$T=\frac{4\pi}{b}$. Hence, we look for periodic solutions of \eqref{Equation}
with the same period $T$, which is equivalent to solving \eqref{Equation} in
the interval $[0,T]$ with the following boundary condition:
\begin{equation}
z\left(  0\right)  -z\left(  T\right)  =0. \label{BC}%
\end{equation}

Denote by $\left\Vert \text{\textperiodcentered}\right\Vert _{C}=\max
_{t\in\lbrack0,T]}\left\vert x\left(  t\right)  \right\vert $ the norm in the
space $C\left(  [0,T],\mathbb{R}\right)  .$

We will establish that the set of solutions of problem
\eqref{Equation}--\eqref{BC}, assuming that it is non-empty, is bounded in the
space $C\left(  [0,T],\mathbb{R}\right)  .$

\begin{theorem}
\label{ThBoundPer}There exists a constant $D>0$ such that every solution
$z($\textperiodcentered$)$ of problem \eqref{Equation}--\eqref{BC} satisfies%
\begin{equation}
\left\Vert z\right\Vert _{C}\leq D. \label{BoundPer}%
\end{equation}

\end{theorem}

\begin{proof}
On the one hand, making use of the equality
\[
F_{z}\left(  t,z\left(  t\right)  \right)  z^{\prime}=\frac{\mathrm{d}%
V}{\mathrm{d}t}(t,z\left(  t\right)  )-\frac{\partial V}{\partial
t}(t,z\left(  t\right)  ),
\]
multiplying \eqref{Equation} by $z^{\prime}$ and integrating over $\left(
0,T\right)  $ we have%
\begin{align}
\int_{0}^{T}\left\vert z^{\prime}\right\vert ^{2}\mathrm{d}t  &  =\int_{0}%
^{T}\left\vert F_{z}(t,z\left(  t\right)  )\right\vert ^{2}\mathrm{d}%
t\nonumber\\
&  =V\left(  T,z\left(  T\right)  \right)  -V\left(  0,z\left(  0\right)
\right)  -\int_{0}^{T}\frac{\partial V}{\partial t}(t,z\left(  t\right)
)\mathrm{d}t\nonumber\\
&  =-\int_{0}^{T}\frac{\partial V}{\partial t}(t,z\left(  t\right)
)\mathrm{d}t, \label{Eq2}%
\end{align}
due to the boundary conditions \eqref{BC}.

On the other hand, since%
\[
\frac{\partial V}{\partial t}(t,z)=\frac{b}{2}F_{0}f\left(  z\right)
\sin\left(  2kz-bt\right)  ,
\]
we have%
\[
F_{z}\left(  t,z\right)  =-\frac{2k}{b}\frac{\partial V}{\partial
t}(t,z)+F_{0}\left(  \cos\left(  kz-\frac{bt}{2}\right)  \right)
^{2}f^{\prime}\left(  z\right)  ,
\]
so that integrating \eqref{Equation} over $(0,T)$ we deduce that
%\begin{align*}
%&  \frac{2k}{b}\int_{0}^{T}\frac{\partial V}{\partial t}(t,z_{\lambda_{n}%
%}\left(  t\right)  )dt\\
%&  =z\left(  0\right)  -z\left(  T\right)  +F_{0}\int_{0}^{T}\left(
%\cos\left(  kz\left(  t\right)  -\frac{bt}{2}\right)  \right)  ^{2}f^{\prime
%}\left(  z\left(  t\right)  \right)  dt\\
%&  =F_{0}\int_{0}^{T}\left(  \cos\left(  kz\left(  t\right)  -\frac{bt}%
%{2}\right)  \right)  ^{2}f^{\prime}\left(  z\left(  t\right)  \right)  dt.
%\end{align*}%
\begin{align*}
&  \frac{2k}{b}\int_{0}^{T}\frac{\partial V}{\partial t}(t,z\left(  t\right)
)\mathrm{d}t\\
&  =z\left(  0\right)  -z\left(  T\right)  +F_{0}\int_{0}^{T}\left(
\cos\left(  kz\left(  t\right)  -\frac{bt}{2}\right)  \right)  ^{2}f^{\prime
}\left(  z\left(  t\right)  \right)  \mathrm{d}t\\
&  =F_{0}\int_{0}^{T}\left(  \cos\left(  kz\left(  t\right)  -\frac{bt}%
{2}\right)  \right)  ^{2}f^{\prime}\left(  z\left(  t\right)  \right)
\mathrm{d}t.
\end{align*}

Substituting this expression in \eqref{Eq2} we have%
\begin{align}
&  \int_{0}^{T}\left\vert F_{z}(t,z\left(  t\right)  )\right\vert
^{2}\mathrm{d}t\nonumber\\
&  =\int_{0}^{T}\left(  -kF_{0}f\left(  z(t)\right)  \sin\left(
2kz(t)-bt\right)  +F_{0}\left(  \cos\left(  kz(t)-\frac{bt}{2}\right)
\right)  ^{2}f^{\prime}\left(  z(t)\right)  \right)  ^{2}\mathrm{d}%
t\nonumber\\
&  =-\frac{bF_{0}}{2k}\int_{0}^{T}\left(  \cos\left(  kz\left(  t\right)
-\frac{bt}{2}\right)  \right)  ^{2}f^{\prime}\left(  z\left(  t\right)
\right)  \mathrm{d}t. \label{EqContrad}%
\end{align}

We state first that that there exists $D_{1}>0$ such that for any solution
$z($\textperiodcentered$)$ of problem \eqref{Equation}--\eqref{BC} there
exists $t_{z}$ for which%
\begin{equation}
\left\vert z(t_{z})\right\vert \leq D_{1}. \label{EstInf}%
\end{equation}

Otherwise, there would exist a sequence of solutions $z_{n}($%
\textperiodcentered$)$ such that%
\begin{equation}
\inf_{t\in\lbrack0,T]}\left\vert z_{n}(t)\right\vert \rightarrow\infty.
\label{InfInfinity}%
\end{equation}
Then integrating equation \eqref{Equation} over $\left(  s,t\right)  $ we
obtain%
\begin{align*}
\left\vert z_{n}(t)-z_{n}(s)\right\vert  &  \leq\int_{s}^{t}\left\vert
F_{z}(r,z_{n}(r))\right\vert \mathrm{d}r\\
&  \leq TF_{0}\left(  k\sup_{r\in\lbrack0,T]}\left\vert f(z_{n}(r))\right\vert
+\sup_{r\in\lbrack0,T]}\left\vert f^{\prime}(z_{n}(r))\right\vert \right)  .
\end{align*}
In view of \eqref{InfInfinity} and the properties of the function $f$ we have%
\begin{equation}
\sup_{t,s\in\lbrack0,T]}\left\vert z_{n}(t)-z_{n}(s)\right\vert \rightarrow
0\text{ as }n\rightarrow\infty. \label{Convergzn}%
\end{equation}

Further, we will analyze each term in equality \eqref{EqContrad} in order to
get a contradiction.

Using the second mean value theorem for integrals it follows the existence of
$\overline{z}_{n}\in\{z_{n}(t),\ t\in\lbrack0,T]\}$ such that%
\[
-\frac{bF_{0}}{2k}\int_{0}^{T}\left(  \cos\left(  kz_{n}\left(  t\right)
-\frac{bt}{2}\right)  \right)  ^{2}f^{\prime}\left(  z_{n}\left(  t\right)
\right)  \mathrm{d}t=-\frac{bF_{0}}{2k}f^{\prime}\left(  \overline{z}%
_{n}\right)  \int_{0}^{T}\left(  \cos\left(  kz_{n}\left(  t\right)
-\frac{bt}{2}\right)  \right)  ^{2}\mathrm{d}t.
\]
Hence, \eqref{EqContrad} yields
\begin{align}
&  \int_{0}^{T}\left(  \cos\left(  kz_{n}\left(  t\right)  -\frac{bt}%
{2}\right)  \right)  ^{2}\mathrm{d}t\label{Ineq1}\\
&  \leq\frac{2k}{bF_{0}}\frac{1}{\left\vert f^{\prime}\left(  \overline{z}%
_{n}\right)  \right\vert }\int_{0}^{T}\left(  -kF_{0}f\left(  z_{n}\left(
t\right)  \right)  \sin\left(  2kz_{n}\left(  t\right)  -bt\right)
+F_{0}\left(  \cos\left(  kz_{n}\left(  t\right)  -\frac{bt}{2}\right)
\right)  ^{2}f^{\prime}\left(  z_{n}\left(  t\right)  \right)  \right)
^{2}\mathrm{d}t\nonumber\\
&  \leq\frac{4kF_{0}}{b}\frac{1}{\left\vert f^{\prime}\left(  \overline{z}%
_{n}\right)  \right\vert }\int_{0}^{T}\left(  k^{2}(f\left(  z_{n}\left(
t\right)  \right)  )^{2}\sin^{2}\left(  2kz_{n}\left(  t\right)  -bt\right)
+\left(  \cos\left(  kz_{n}\left(  t\right)  -\frac{bt}{2}\right)  \right)
^{4}(f^{\prime}\left(  z_{n}\left(  t\right)  \right)  )^{2}\right)
\mathrm{d}t\nonumber\\
&  \leq\frac{4kF_{0}}{b}\frac{1}{\left\vert f^{\prime}\left(  \overline{z}%
_{n}\right)  \right\vert }\int_{0}^{T}\left(  k^{2}(f\left(  z_{n}\left(
t\right)  \right)  )^{2}+(f^{\prime}\left(  z_{n}\left(  t\right)  \right)
)^{2}\right)  \mathrm{d}t\nonumber\\
&  =\frac{4kF_{0}}{b}\frac{1}{\left\vert f^{\prime}\left(  \overline{z}%
_{n}\right)  \right\vert }\left(  k^{2}(f\left(  \widetilde{z}_{n}\right)
)^{2}+(f^{\prime}\left(  \widetilde{z}_{n}\right)  )^{2}\right)  ,\nonumber
\end{align}
for some $\widetilde{z}_{n}\in\{z_{n}(t),\ t\in\lbrack0,T]\}$. We observe that
$f^{\prime}\left(  \overline{z}_{n}\right)  \not =0$, because otherwise from
\eqref{EqContrad} we would obtain that%
\[
\int_{0}^{T}\left\vert F_{z}(t,z_{n}\left(  t\right)  )\right\vert
^{2}\mathrm{d}t=\int_{0}^{T}\left\vert z_{n}^{\prime}\left(  t\right)
\right\vert ^{2}\mathrm{d}t=0,
\]
so $z_{n}^{\prime}\left(  t\right)  \equiv0$ and $z_{n}($\textperiodcentered
$)$ would be a fixed point, which is impossible by assumption.

The rigth-hand side of \eqref{Ineq1} converges to $0$ in light of assumption
\eqref{Cond4Functf} as \eqref{Convergzn} implies that $\frac{\widetilde{z}%
_{n}}{\overline{z}_{n}}\rightarrow1$. But we can prove that
\[
\int_{0}^{T}\left(  \cos\left(  kz_{n}\left(  t\right)  -\frac{bt}{2}\right)
\right)  ^{2}\mathrm{d}t\rightarrow\frac{T}{2},
\]
obtaining in this way the desired contradiction. Indeed, it is clear that%
\[
\int_{0}^{T}\left(  \cos\left(  kz_{n}\left(  0\right)  -\frac{bt}{2}\right)
\right)  ^{2}\mathrm{d}t=\frac{T}{2},
\]
and \eqref{Convergzn} implies that%
\[
\left\vert \left(  \cos\left(  kz_{n}\left(  t\right)  -\frac{bt}{2}\right)
\right)  ^{2}-\left(  \cos\left(  kz_{n}\left(  0\right)  -\frac{bt}%
{2}\right)  \right)  ^{2}\right\vert \rightarrow0\text{ uniformly in }[0,T],
\]

so%

\begin{align*}
&  \left\vert \int_{0}^{T}\left(  \cos\left(  kz_{n}\left(  t\right)
-\frac{bt}{2}\right)  \right)  ^{2}dt-\frac{T}{2}\right\vert \\
&  =\left\vert \int_{0}^{T}\left(  \left(  \cos\left(  kz_{n}\left(  t\right)
-\frac{bt}{2}\right)  \right)  ^{2}-\left(  \cos\left(  kz_{n}\left(
0\right)  -\frac{bt}{2}\right)  \right)  ^{2}\right)  dt\right\vert
\rightarrow0.
\end{align*}

Once we have proved that \eqref{EstInf} holds true, the statement of the
theorem is deduced easily. Integrating equation \eqref{Equation} over the
interval $(t_{z},t)$ we get%
\begin{align*}
\left\vert z(t)\right\vert  &  \leq\left\vert z(t_{z})\right\vert +\int%
_{t_{z}}^{t}\left\vert F_{z}(s,z(s))\right\vert \mathrm{d}s\\
&  \leq D_{1}+F_{0}T\left(  k\sup_{r\in\mathbb{R}}\left\vert f(r)\right\vert
+\sup_{r\in\mathbb{R}}\left\vert f^{\prime}(r)\right\vert \right)  =D_{2},
\end{align*}
for any $t\in\lbrack t_{z},T]$. Finally, since $z(0)=z(T)$, for $t\in
\lbrack0,t_{z}]$ we get%
\begin{align*}
\left\vert z(t)\right\vert  &  \leq\left\vert z(0)\right\vert +\int_{0}%
^{t}\left\vert F_{z}(s,z(s))\right\vert \mathrm{d}s\\
&  \leq D_{2}+F_{0}T\left(  k\sup_{r\in\mathbb{R}}\left\vert f(r)\right\vert
+\sup_{r\in\mathbb{R}}\left\vert f^{\prime}(r)\right\vert \right)  =D,
\end{align*}
which concludes the proof.
\end{proof}

\bigskip

\begin{corollary}
Every periodic solution with period $T$ belongs to the ball of radius $D$
centered at $0$ of the space of continuous bounded functions $C_{\mathrm{b}%
}(\mathbb{R})$, whose norm is $\left\Vert u\right\Vert _{C_{\mathrm{b}}}%
=\sup_{t\in\mathbb{R}}\left\vert u(t)\right\vert .$
\end{corollary}

The following criterion is useful in order to check that the equation has no
fixed points.

\begin{lemma}
\label{FixedCriterion}$\overline{z}$ is a fixed point of \eqref{Equation} if
and only if%
\begin{equation}
f(\overline{z})=f^{\prime}(\overline{z})=0. \label{Eqffprima}%
\end{equation}

\end{lemma}

\begin{proof}
Let $\overline{z}$ be a fixed point, so $F_{z}(s,\overline{z})=0$ for any $s$,
that is,%
\[
kf\left(  \overline{z}\right)  \sin\left(  2k\overline{z}-bt\right)  =\left(
\cos\left(  k\overline{z}-\frac{bt}{2}\right)  \right)  ^{2}f^{\prime}\left(
\overline{z}\right)  \ \forall t.
\]
Assume that \eqref{Eqffprima} is not true. If $f\left(  \overline{z}\right)
=0,\ f^{\prime}\left(  \overline{z}\right)  \not =0$ (respectively, $f\left(
\overline{z}\right)  \not =0,\ f^{\prime}\left(  \overline{z}\right)  =0$),
then $\cos\left(  k\overline{z}-\frac{bt}{2}\right)  =0\ $for any $t$
(respectively, $\sin\left(  2k\overline{z}-bt\right)  =0$), which is not
possible. If $f\left(  \overline{z}\right)  \not =0,\ f^{\prime}\left(
\overline{z}\right)  \not =0$, then $\tan\left(  k\overline{z}-\frac{b}%
{2}t\right)  =\frac{f^{\prime}\left(  \overline{z}\right)  }{2kf\left(
\overline{z}\right)  }\ $for any $t$, which, again, cannot occur. Therefore,
if $\overline{z}$ is a fixed point, then \eqref{Eqffprima} holds.

The converse statement is straightforward.
\end{proof}

\bigskip

\subsection{Existence of periodic solutions}

We will prove the existence of periodic solutions in the aforementioned
particular cases.

\subsubsection{Case 1: $f(z)=\frac{z_{0}^{2}}{z_{0}^{2}+z^{2}}$}

We shall prove in this section that problem \eqref{Equation}, \eqref{BC} is
solvable in the particular case where $f(z)=\frac{z_{0}^{2}}{z_{0}^{2}+z^{2}}%
$, where $z_{0}>0.$

First of all, let us define the functions $p{\colon}C([0,T],\mathbb{R})\times
C([0,T],\mathbb{R})\rightarrow$$C([0,T],\mathbb{R})$, $l{\colon}%
C([0,T],\mathbb{R})\times C([0,T],\mathbb{R})\rightarrow\mathbb{R}$ by%
\begin{align}
p(z,y)(t)  &  =-y\left(  t\right)  ,\ \forall t\in\lbrack0,T],\nonumber\\
l(z,y)  &  =y(0)-y(T). \label{Pairpl}%
\end{align}
It is obvious that the maps $p(z,$\textperiodcentered$){\colon}$%
$C([0,T],\mathbb{R})\rightarrow C([0,T],\mathbb{R})$, $l(z,$%
\textperiodcentered$){\colon}$$C([0,T],\mathbb{R})\rightarrow\mathbb{R}$ are
linear for any fixed $z\in$$C([0,T],\mathbb{R})$. Also, for any $y,z\in
$$C([0,T],\mathbb{R})$ and $t$$\in\lbrack0,T]$ we have
\begin{align*}
\left\vert p(z,y)(t)\right\vert  &  \leq\left\vert y\left(  t\right)
\right\vert ,\\
\left\vert l(z,y)\right\vert  &  \leq2\left\Vert y\right\Vert _{C}.
\end{align*}

Consider the following problem%
\begin{equation}
\frac{\mathrm{d}y}{\mathrm{d}t}=p(z,y)(t)+q(t),\ l(z,y)=c_{0},
\label{Problempq}%
\end{equation}
where $q$$\in C([0,T],\mathbb{R})$ and $c_{0}\in\mathbb{R}.$ Every solution of
this problem satisfies%
\begin{equation}
y(t)=\exp\left(  -t\right)  y(0)+\int_{0}^{t}\exp\left(  -\left(  t-s\right)
\right)  q(s)\mathrm{d}s. \label{SolLinear}%
\end{equation}
Using the boundary conditions we have%
\[
y\left(  T\right)  =\exp\left(  -T\right)  (y\left(  T\right)  +c_{0}%
)+\int_{0}^{T}\exp\left(  -(T-s)\right)  q(s)\mathrm{d}s,
\]
so%
\[
\left\vert y\left(  T\right)  \right\vert \leq\frac{1}{1-\exp\left(
-T\right)  }\left(  \left\vert c_{0}\right\vert +\int_{0}^{T}\left\vert
q\left(  s\right)  \right\vert \mathrm{d}s\right)
\]
and%
\[
\left\vert y\left(  0\right)  \right\vert \leq\left\vert y\left(  T\right)
\right\vert +\left\vert c_{0}\right\vert \leq\frac{1}{1-\exp\left(  -T\right)
}\left(  \left\vert c_{0}\right\vert +\int_{0}^{T}\left\vert q\left(
s\right)  \right\vert \mathrm{d}s\right)  +\left\vert c_{0}\right\vert .
\]
Therefore, by \eqref{SolLinear} there exists $\beta>0$ such that%
\begin{equation}
\left\Vert y\right\Vert _{C}\leq\beta(\left\vert c_{0}\right\vert +\left\Vert
q\right\Vert _{L^{1}(0,T;\mathbb{R})}). \label{Est1}%
\end{equation}

Further, we study the boundary-value problem%
\begin{equation}
\left\{
\begin{array}
[c]{c}%
\dfrac{\mathrm{d}z}{\mathrm{d}t}=p(z,z)(t)+\lambda\left(  F_{z}(t,z\left(
t\right)  )-p\left(  z,z\right)  \left(  t\right)  \right)  ,\\
z(0)-z\left(  T\right)  =0,
\end{array}
\right.  \label{BVProblemAux}%
\end{equation}
for $\lambda\in\left(  0,1\right)  $.

\begin{theorem}
\label{ThCotaZetaLambda}There exists $\rho>0$ such that for any $\lambda
\in\left(  0,1\right)  $ and any solution $z_{\lambda}\left(
\text{\textperiodcentered}\right)  $ to problem \eqref{BVProblemAux} it holds
that%
\begin{equation}
\left\Vert z_{\lambda}\right\Vert _{C}\leq\rho. \label{CotaZLambda}%
\end{equation}

\end{theorem}

\begin{proof}
First, we shall prove the existence of $\rho^{\prime}>0$ such that for any
$\lambda\in\left(  0,1\right)  $ and any solution $z_{\lambda}\left(
\text{\textperiodcentered}\right)  $ to problem \eqref{BVProblemAux} the
estimate
\begin{equation}
\min_{t\in\lbrack0,T]}\left\vert z_{\lambda}\left(  t\right)  \right\vert
<\rho^{\prime} \label{EstRoPrima}%
\end{equation}
is satisfied. By contradiction, if this is not true, then there exists a
sequence of solutions $z_{\lambda_{n}}\left(  \text{\textperiodcentered
}\right)  $ such that%
\[
\min_{t\in\lbrack0,T]}\left\vert z_{\lambda_{n}}\left(  t\right)  \right\vert
\rightarrow\infty.
\]
Without loss of generality we can assume that $\min_{t\in\lbrack
0,T]}z_{\lambda_{n}}\left(  t\right)  \rightarrow+\infty$.

Multiplying the equation in \eqref{BVProblemAux} by $z_{0}^{2}+z_{\lambda_{n}%
}^{2}(t)$ and integrating over $(0,T)$ we have%
\begin{align*}
&  -\left(  1-\lambda_{n}\right)  \int_{0}^{T}p(z_{\lambda_{n}},z_{\lambda
_{n}})(t)(z_{0}^{2}+z_{\lambda_{n}}^{2}(t))\mathrm{d}t\\
&  =-kF_{0}z_{0}^{2}\lambda_{n}\int_{0}^{T}\sin\left(  2kz_{\lambda_{n}%
}\left(  t\right)  -bt\right)  \mathrm{d}t+\lambda_{n}F_{0}\int_{0}^{T}\left(
\cos\left(  kz_{\lambda_{n}}\left(  t\right)  -\frac{bt}{2}\right)  \right)
^{2}f^{\prime}\left(  z_{\lambda_{n}}\left(  t\right)  \right)  (z_{0}%
^{2}+z_{\lambda_{n}}^{2}(t))\mathrm{d}t.
\end{align*}

There exist $z_{\lambda_{n}}^{\ast}=z_{\lambda_{n}}(t_{n}^{\ast}%
),\ \overline{z}_{\lambda_{n}}=z_{\lambda_{n}}(\overline{t}_{n})$ such that
\begin{align*}
&  T(1-\lambda_{n})z_{\lambda_{n}}^{\ast}(z_{0}^{2}+(z_{\lambda_{n}}^{\ast
})^{2})\\
&  =-kF_{0}\lambda_{n}z_{0}^{2}\int_{0}^{T}\sin\left(  2kz_{\lambda_{n}%
}\left(  t\right)  -bt\right)  \mathrm{d}t+\lambda_{n}F_{0}f^{\prime}\left(
\overline{z}_{\lambda_{n}}\right)  (z_{0}^{2}+(\overline{z}_{\lambda_{n}}%
)^{2})\int_{0}^{T}\left(  \cos\left(  kz_{\lambda_{n}}\left(  t\right)
-\frac{bt}{2}\right)  \right)  ^{2}\mathrm{d}t\\
&  =-kF_{0}\lambda_{n}z_{0}^{2}\int_{0}^{T}\sin\left(  2kz_{\lambda_{n}%
}\left(  t\right)  -bt\right)  \mathrm{d}t-2\lambda_{n}F_{0}z_{0}^{2}%
\frac{\overline{z}_{\lambda_{n}}}{z_{0}^{2}+(\overline{z}_{\lambda_{n}})^{2}%
}\int_{0}^{T}\left(  \cos\left(  kz_{\lambda_{n}}\left(  t\right)  -\frac
{bt}{2}\right)  \right)  ^{2}\mathrm{d}t,
\end{align*}
so%
\begin{align}
&  T(1-\lambda_{n})\frac{z_{\lambda_{n}}^{\ast}}{\overline{z}_{\lambda_{n}}%
}(z_{0}^{2}+(z_{\lambda_{n}}^{\ast})^{2})(z_{0}^{2}+(\overline{z}_{\lambda
_{n}})^{2})\label{EqContradiction}\\
&  =-kF_{0}\lambda_{n}z_{0}^{2}\frac{z_{0}^{2}+(\overline{z}_{\lambda_{n}%
})^{2}}{\overline{z}_{\lambda_{n}}}\int_{0}^{T}\sin\left(  2kz_{\lambda_{n}%
}\left(  t\right)  -bt\right)  \mathrm{d}t-2\lambda_{n}F_{0}z_{0}^{2}\int%
_{0}^{T}\left(  \cos\left(  kz_{\lambda_{n}}\left(  t\right)  -\frac{bt}%
{2}\right)  \right)  ^{2}\mathrm{d}t.\nonumber
\end{align}

Integrating the equation in \eqref{BVProblemAux} over $\left(  0,T\right)  $
and using $z_{\lambda}\left(  T\right)  =z_{\lambda}\left(  0\right)  $ we
obtain for any $\lambda$ that
\begin{equation}
-(1-\lambda)\int_{0}^{T}p(z_{\lambda},z_{\lambda})(s)\mathrm{d}s=\lambda
\int_{0}^{T}F_{z}(s,z_{\lambda}\left(  s\right)  )\mathrm{d}s. \label{EqPZ}%
\end{equation}
Since $z_{\lambda_{n}}\left(  t\right)  >0$ and $p(z_{\lambda_{n}}%
,z_{\lambda_{n}})(s)=-z_{\lambda_{n}}(s)$, we have%
\begin{align}
-(1-\lambda_{n})\int_{0}^{T}p(z_{\lambda_{n}},z_{\lambda_{n}})(s)\mathrm{d}s
&  =(1-\lambda_{n})\int_{0}^{T}z_{\lambda_{n}}(s)\mathrm{d}s\nonumber\\
&  =(1-\lambda_{n})\int_{0}^{T}\left\vert p(z_{\lambda_{n}},z_{\lambda_{n}%
})(s)\right\vert \mathrm{d}s=\lambda_{n}\int_{0}^{T}F_{z}(s,z_{\lambda_{n}%
}\left(  s\right)  )\mathrm{d}s. \label{EqIntp}%
\end{align}

Therefore, we deduce from \eqref{BVProblemAux} and \eqref{EqIntp} that%
\begin{align}
\left\vert z_{\lambda_{n}}(t)-z_{\lambda_{n}}(0)\right\vert  &  \leq\int%
_{0}^{T}\left\vert (1-\lambda_{n})p(z_{\lambda_{n}},z_{\lambda_{n}%
})(s)+\lambda_{n}F_{z}(s,z_{\lambda_{n}}(s))\right\vert \mathrm{d}s\nonumber\\
&  \leq2\int_{0}^{T}\left\vert F_{z}(s,z_{\lambda_{n}}(s))\right\vert
\mathrm{d}s\nonumber\\
&  \leq2\int_{0}^{T}\left(  kF_{0}\left\vert f\left(  z_{\lambda_{n}%
}(s)\right)  \right\vert +F_{0}\left\vert f^{\prime}\left(  z_{\lambda_{n}%
}(s)\right)  \right\vert \right)  \mathrm{d}s\nonumber\\
&  =2F_{0}T(k\left\vert f\left(  \widetilde{z}_{\lambda_{n}}\right)
\right\vert +\left\vert f^{\prime}\left(  \widetilde{z}_{\lambda_{n}}\right)
\right\vert )\nonumber\\
&  =2F_{0}Tz_{0}^{2}\left(  \frac{k}{z_{0}^{2}+(\widetilde{z}_{\lambda_{n}%
})^{2}}+\frac{2\widetilde{z}_{\lambda_{n}}}{\left(  z_{0}^{2}+(\widetilde{z}%
_{\lambda_{n}})^{2}\right)  ^{2}}\right)  ,\text{ }\forall t\in\lbrack0,T],
\label{EstDiffZ}%
\end{align}
for some $\widetilde{z}_{\lambda_{n}}=z_{\lambda_{n}}(\widetilde{t}_{n})$.
This implies, in particular, that%
\begin{equation}
\sup_{t,s\in\lbrack0,T]}\left\vert z_{\lambda_{n}}(t)-z_{\lambda_{n}%
}(s)\right\vert \rightarrow0, \label{Convergz}%
\end{equation}
as $n\rightarrow\infty$.

If we proved that%
\begin{equation}
\frac{z_{0}^{2}+(\overline{z}_{\lambda_{n}})^{2}}{\overline{z}_{\lambda_{n}}%
}\int_{0}^{T}\sin\left(  2kz_{\lambda_{n}}\left(  t\right)  -bt\right)
\mathrm{d}t\rightarrow0, \label{Converg1}%
\end{equation}%
\begin{equation}
\int_{0}^{T}\left(  \cos\left(  kz_{\lambda_{n}}\left(  t\right)  -\frac
{bt}{2}\right)  \right)  ^{2}\mathrm{d}t\rightarrow\frac{T}{2},
\label{Converg2}%
\end{equation}
then \eqref{EqContradiction} would imply that%
\[
\frac{1-\lambda_{n}}{\lambda_{n}}\frac{z_{\lambda_{n}}^{\ast}}{\overline
{z}_{\lambda_{n}}}(z_{0}^{2}+(z_{\lambda_{n}}^{\ast})^{2})(z_{0}%
^{2}+(\overline{z}_{\lambda_{n}})^{2})\rightarrow-F_{0}z_{0}^{2},
\]
which is not possible as $\frac{1-\lambda_{n}}{\lambda_{n}}\frac
{z_{\lambda_{n}}^{\ast}}{\overline{z}_{\lambda_{n}}}(z_{0}^{2}+(z_{\lambda
_{n}}^{\ast})^{2})(z_{0}^{2}+(\overline{z}_{\lambda_{n}})^{2})>0$ for all $n$.
Thus, \eqref{EstRoPrima} would be true.

In order to check \eqref{Converg1}, making use of the fact that
\[
\int_{0}^{T}\sin\left(  2kz_{\lambda_{n}}\left(  0\right)  -bt\right)
\mathrm{d}t=0,
\]
we obtain by \eqref{EstDiffZ} that%
\begin{align*}
&  \frac{z_{0}^{2}+(\overline{z}_{\lambda_{n}})^{2}}{\overline{z}_{\lambda
_{n}}}\left\vert \int_{0}^{T}\sin\left(  2kz_{\lambda_{n}}\left(  t\right)
-bt\right)  \mathrm{d}t\right\vert \\
&  \leq\frac{z_{0}^{2}+(\overline{z}_{\lambda_{n}})^{2}}{\overline{z}%
_{\lambda_{n}}}\int_{0}^{T}\left\vert \sin\left(  2kz_{\lambda_{n}}\left(
t\right)  -bt\right)  -\sin\left(  2kz_{\lambda_{n}}\left(  0\right)
-bt\right)  \right\vert \mathrm{d}t\\
&  \leq\frac{z_{0}^{2}+(\overline{z}_{\lambda_{n}})^{2}}{\overline{z}%
_{\lambda_{n}}}\int_{0}^{T}\left\vert \cos(\beta_{n}(t))\right\vert
2k\left\vert z_{\lambda_{n}}\left(  t\right)  -z_{\lambda_{n}}\left(
0\right)  \right\vert \mathrm{d}t\\
&  \leq\frac{z_{0}^{2}+(\overline{z}_{\lambda_{n}})^{2}}{\overline{z}%
_{\lambda_{n}}}4kT^{2}F_{0}z_{0}^{2}\left(  \frac{k}{z_{0}^{2}+(\widetilde{z}%
_{\lambda_{n}})^{2}}+\frac{2\widetilde{z}_{\lambda_{n}}}{\left(  z_{0}%
^{2}+(\widetilde{z}_{\lambda_{n}})^{2}\right)  ^{2}}\right)  \rightarrow0,
\end{align*}
as from \eqref{Convergz} we can easily see that%
\begin{align*}
\frac{z_{0}^{2}+(\overline{z}_{\lambda_{n}})^{2}}{z_{0}^{2}+(\widetilde{z}%
_{\lambda_{n}})^{2}}  &  \rightarrow1,\\
\frac{\widetilde{z}_{\lambda_{n}}}{\overline{z}_{\lambda_{n}}}  &
\rightarrow1.
\end{align*}
In order to obtain \eqref{Converg2} it suffices to see that%
\[
\int_{0}^{T}\left(  \cos\left(  kz_{\lambda_{n}}\left(  0\right)  -\frac
{bt}{2}\right)  \right)  ^{2}\mathrm{d}t=\frac{T}{2}%
\]
and that \eqref{Convergz} implies%
\[
\left\vert \left(  \cos\left(  kz_{\lambda_{n}}\left(  t\right)  -\frac{bt}%
{2}\right)  \right)  ^{2}-\left(  \cos\left(  kz_{\lambda_{n}}\left(
0\right)  -\frac{bt}{2}\right)  \right)  ^{2}\right\vert \rightarrow0\text{
uniformly in }[0,T],
\]
so%
\begin{align*}
&  \left\vert \int_{0}^{T}\left(  \cos\left(  kz_{\lambda_{n}}\left(
t\right)  -\frac{bt}{2}\right)  \right)  ^{2}\mathrm{d}t-\frac{T}%
{2}\right\vert \\
&  =\left\vert \int_{0}^{T}\left(  \left(  \cos\left(  kz_{\lambda_{n}}\left(
t\right)  -\frac{bt}{2}\right)  \right)  ^{2}-\left(  \cos\left(
kz_{\lambda_{n}}\left(  0\right)  -\frac{bt}{2}\right)  \right)  ^{2}\right)
\mathrm{d}t\right\vert \rightarrow0.
\end{align*}

Finally, let us prove that \eqref{EstRoPrima} implies \eqref{CotaZLambda}.

It follows from \eqref{EstRoPrima}\ that for every solution $z_{\lambda}%
($\textperiodcentered$)$ of problem \eqref{BVProblemAux} there exists a moment
of time $t_{z_{\lambda}}$ such that $\left\vert z_{\lambda}(t_{z_{\lambda}%
})\right\vert <\rho^{\prime}$. Integrating the equation in
\eqref{BVProblemAux} over $\left(  t_{z_{\lambda}},t\right)  $ we have%
\[
z_{\lambda}(t)=z_{\lambda}(t_{z_{\lambda}})+\int_{t_{z_{\lambda}}}^{t}\left(
-\left(  1-\lambda\right)  z_{\lambda}(s)+\lambda F_{z}(t,z_{\lambda}\left(
s\right)  )\right)  \mathrm{d}s.
\]
Hence,%
\begin{align*}
\left\vert z_{\lambda}(t)\right\vert  &  \leq\rho^{\prime}+\int_{t_{z_{\lambda
}}}^{t}\left(  \left(  1-\lambda\right)  \left\vert z_{\lambda}(s)\right\vert
+\lambda\left\vert F_{z}(t,z_{\lambda}\left(  s\right)  )\right\vert \right)
\mathrm{d}s\\
&  \leq\rho^{\prime}+\int_{t_{z_{\lambda}}}^{t}\left(  \left\vert z_{\lambda
}(s)\right\vert +kF_{0}\max_{u\in\mathbb{R}}\left\vert f\left(  u\right)
\right\vert +F_{0}\max_{u\in\mathbb{R}}\left\vert f^{\prime}\left(  u\right)
\right\vert \right)  \mathrm{d}s\\
&  \leq\rho^{\prime}+T\left(  kF_{0}\max_{u\in\mathbb{R}}\left\vert f\left(
u\right)  \right\vert +F_{0}\max_{u\in\mathbb{R}}\left\vert f^{\prime}\left(
u\right)  \right\vert \right)  +\int_{t_{z_{\lambda}}}^{t}\left\vert
z_{\lambda}(s)\right\vert \mathrm{d}s\\
&  =R+\int_{t_{z_{\lambda}}}^{t}\left\vert z_{\lambda}(s)\right\vert
\mathrm{d}s.
\end{align*}
Applying Gronwall's lemma we get%
\[
\left\vert z_{\lambda}(t)\right\vert \leq Re^{t-t_{z_{\lambda}}}\leq
Re^{T}\text{ for all }t\in\lbrack t_{z_{\lambda}},T]\text{.}%
\]
Noting that $z_{\lambda}(T)=z_{\lambda}(0)$ implies%
\[
\left\vert z_{\lambda}(0)\right\vert =\left\vert z_{\lambda}(T)\right\vert
\leq Re^{T},
\]
we can repeat the same argument in order to show the existence of $\rho\geq
Re^{T}$ such that%
\[
\left\vert z_{\lambda}(t)\right\vert \leq\rho\text{ for all }t\in
\lbrack0,t_{z_{\lambda}}]\text{.}%
\]
The proof of \eqref{CotaZLambda} is now complete.
\end{proof}

\bigskip

Let us recall a general result on existence of solutions for boundary-value
problems, which was proved in \cite{Kiguradze97}.

Let $I=[a,b]\subset\mathbb{R}$, $g{\colon}C(I,\mathbb{R}^{n})\rightarrow
L^{1}(I,\mathbb{R}^{n})$, $h{\colon}C(I,\mathbb{R}^{n})\rightarrow
\mathbb{R}^{n}$ be continuous operators such that for any $\rho>0$ one has%
\begin{equation}
\sup\{\left\Vert g(x)(\text{\textperiodcentered})\right\Vert _{\mathbb{R}^{n}%
}{\colon}x\in C(I,\mathbb{R}^{n}),\ \left\Vert x\right\Vert _{C}\leq\rho\}\in
L^{1}(I,\mathbb{R}^{n}), \label{supf}%
\end{equation}%
\begin{equation}
\sup\{\left\Vert h(x)\right\Vert _{\mathbb{R}^{n}}{\colon}x\in C(I,\mathbb{R}%
^{n}),\ \left\Vert x\right\Vert _{C}\leq\rho\}<\infty, \label{suph}%
\end{equation}
where, as before, $\left\Vert x\right\Vert _{C}=\max_{t\in I}\left\Vert
x(t)\right\Vert _{\mathbb{R}^{n}}$, and consider the boundary-value problem%
\begin{equation}
\left\{
\begin{array}
[c]{c}%
\dfrac{\mathrm{d}x}{\mathrm{d}t}=g(x)(t),\\
h(x)=0.
\end{array}
\right.  \label{BVProblemf}%
\end{equation}

We also consider a pair $(p,l)$ of continuous operators $p{\colon
}C(I,\mathbb{R}^{n})\times C(I,\mathbb{R}^{n})\rightarrow L^{1}(I,\mathbb{R}%
^{n})$, $l{\colon}C(I,\mathbb{R}^{n})\times C(I,\mathbb{R}^{n})\rightarrow
\mathbb{R}^{n}$ satisfying:

\begin{enumerate}
\item for any fixed $x\in C(I,\mathbb{R}^{n})$ the operators $p(x,$%
\textperiodcentered$){\colon}C(I,\mathbb{R}^{n})\rightarrow L^{1}%
(I,\mathbb{R}^{n})$, $l(x,$\textperiodcentered$){\colon}C(I,\mathbb{R}%
^{n})\rightarrow\mathbb{R}^{n}$ are linear;

\item for any $x,y\in C(I,\mathbb{R}^{n})$ we have%
\[
\left\Vert p(x,y)(t)\right\Vert _{\mathbb{R}^{n}}\leq\alpha(t,\left\Vert
x\right\Vert _{C})\left\Vert y\right\Vert _{C},
\]%
\[
\left\Vert l(x,y)(t)\right\Vert _{\mathbb{R}^{n}}\leq\alpha_{0}(\left\Vert
x\right\Vert _{C})\left\Vert y\right\Vert _{C},
\]
with $\alpha_{0}:\mathbb{R}^{+}\rightarrow\mathbb{R}^{+}$ being nondecreasing
and $\alpha:\mathbb{R}^{+}\rightarrow\mathbb{R}^{+}$ being nondecreasing in
the second argument and integrable with respect to the first one;

\item there exists $\beta>0$ such that for any $x\in C(I,\mathbb{R}^{n})$,
$q\in L^{1}(I,\mathbb{R}^{n})$ and $c_{0}\in\mathbb{R}^{n}$ every solution of
the boundary-value problem%
\[
\frac{\mathrm{d}y}{\mathrm{d}t}=p(x,y)(t)+q(t),\ l(x,y)=c_{0},
\]
satisfies the estimate%
\[
\left\Vert y\right\Vert _{C}\leq\beta(\left\Vert c_{0}\right\Vert
_{\mathbb{R}^{n}}+\left\Vert q\right\Vert _{L^{1}(I,\mathbb{R}^{n})}).
\]

\end{enumerate}

\begin{theorem}
\label{ThKiguradze}\cite[Theorem 1]{Kiguradze97} Suppose that there exists
$\rho>0$ such that for any $\lambda\in(0,1)$ the boundary-value problem%
\[
\left\{
\begin{array}
[c]{c}%
\dfrac{\mathrm{d}x}{\mathrm{d}t}=p(x,x)(t)+\lambda(g(x)(t)-p(x,x)(t)),\ \\
l(x,x)=\lambda(l(x,x)-h(x)),
\end{array}
\right.
\]
satisfies the estimate%
\begin{equation}
\left\Vert x\right\Vert _{C}\leq\rho. \label{EstRo}%
\end{equation}

Then problem \eqref{BVProblemf} possesses at least one solution.
\end{theorem}

Applying this theorem we can prove the existence of periodic solutions for
equation \eqref{Equation}.

\begin{theorem}
\label{ThPerPolynomial}The boundary-value problem \eqref{Equation}, \eqref{BC}
possesses at least one solution, which can be extended to a periodic solution
of equation \eqref{Equation}.
\end{theorem}

\begin{proof}
We put $I=[0,T]$, $n=1$, $x=z$, $h(z)=z(0)-z(T)$ and $g(z)(t)=F_{z}(t,z(t))$.
It is not difficult to check that \eqref{supf}-\eqref{suph} hold true. The
pair $(p,l)$ was defined in \eqref{Pairpl} and the properties 1--3 were proved
to be true above. Also, by Theorem \ref{ThCotaZetaLambda} we know that
\eqref{EstRo} is satisfied. Therefore, Theorem \ref{ThKiguradze} proves the
solvability of problem \eqref{Equation}, \eqref{BC}.

Finally, we can define the solution $z($\textperiodcentered$)$ on $[0,\infty)$
by putting $z(t+nT)=z(t)$, for any $t\in\lbrack0,T]$ and $n\in\mathbb{N}$.
This function is a solution of \eqref{Equation} as well due to the fact that
the function $t\mapsto F_{z}(t,z)$ is periodic with period $T$ for any fixed
$z$.
\end{proof}

\bigskip

Using Theorem \ref{ThBoundPer} it is easy to see that the set of all periodic
solutions of equation \eqref{Equation} is bounded.

\begin{theorem}
There exists a constant $D>0$ such that every solution $z($\textperiodcentered
$)$ of problem \eqref{Equation}, \eqref{BC} satisfies%
\[
\left\Vert z\right\Vert _{C}\leq D.
\]

\end{theorem}

\begin{proof}
It suffices to verify that the assumptions of Theorem \ref{ThBoundPer} are
satisfied for $f(z)=\frac{z_{0}^{2}}{z_{0}^{2}+z^{2}}$.

Conditions 1--3 are obvious, and condition 5 follows from Lemma
\ref{FixedCriterion} as $f(z)\not =0$ for all $z.$

It remains to check the fourth one. Let $z_{n},v_{n}\rightarrow\infty$,
$\frac{v_{n}}{z_{n}}\rightarrow1$, as $n\rightarrow\infty$, and $f^{\prime
}\left(  z_{n}\right)  \not =0$. Then%
\[
\frac{f^{2}(v_{n})}{f^{\prime}(z_{n})}=\frac{z_{0}^{4}}{(z_{0}^{2}+v_{n}%
^{2})^{2}}\frac{(z_{0}^{2}+z_{n}^{2})^{2}}{z_{0}^{2}(-2z_{n})}\rightarrow0,
\]%
\[
\frac{(f^{\prime}(v_{n}))^{2}}{f^{\prime}(z_{n})}=\frac{4z_{0}^{4}v_{n}^{2}%
}{(z_{0}^{2}+v_{n}^{2})^{4}}\frac{(z_{0}^{2}+z_{n}^{2})^{2}}{z_{0}^{2}%
(-2z_{n})}\rightarrow0.
\]

\end{proof}

\bigskip

\begin{corollary}
Every periodic solution with period $T$ belongs to the ball of radius $D$
centered at $0$ of the space of continuous bounded functions $C_{\mathrm{b}%
}(\mathbb{R}).$
\end{corollary}

\subsubsection{Case 2: $f(z)=\exp\left(  -2\frac{z^{2}}{z_{0}^{2}}\right)  $}

We will consider in this section the function $f(z)=\exp\left(  -2\frac{z^{2}%
}{z_{0}^{2}}\right)  $, where $z_{0}>0.$

Let us prove that Theorem \ref{ThCotaZetaLambda} is true in this case as well.

\begin{theorem}
\label{ThCotaZetaLambda2}There exists $\rho>0$ such that for any $\lambda
\in\left(  0,1\right)  $ and any solution $z_{\lambda}\left(
\text{\textperiodcentered}\right)  $ to problem \eqref{BVProblemAux} it holds
that%
\begin{equation}
\left\Vert z_{\lambda}\right\Vert _{C}\leq\rho. \label{CotaZLambda2}%
\end{equation}

\end{theorem}

\begin{proof}
First, we need to check the existence of $\rho^{\prime}>0$ such that for any
$\lambda\in\left(  0,1\right)  $ and any solution $z_{\lambda}\left(
\text{\textperiodcentered}\right)  $ to problem \eqref{BVProblemAux} the
estimate
\begin{equation}
\min_{t\in\lbrack0,T]}\left\vert z_{\lambda}\left(  t\right)  \right\vert
<\rho^{\prime} \label{EstRoPrima2}%
\end{equation}
is satisfied. If \eqref{EstRoPrima2} is not true, then there exists a sequence
of solutions $z_{\lambda_{n}}\left(  \text{\textperiodcentered}\right)  $ such
that%
\[
\min_{t\in\lbrack0,T]}\left\vert z_{\lambda_{n}}\left(  t\right)  \right\vert
\rightarrow\infty.
\]
Without loss of generality we can assume that $\min_{t\in\lbrack
0,T]}z_{\lambda_{n}}\left(  t\right)  \rightarrow+\infty$.

We use the same functions $p,l$ as in the previous case. Multiplying the
equation in \eqref{BVProblemAux} by $z_{\lambda_{n}}(t)\exp(2z_{\lambda_{n}%
}^{2}(t)/z_{0}^{2})$ and integrating over $(0,T)$ we have%
\begin{align*}
&  \int_{0}^{T}z_{\lambda_{n}}^{\prime}(t)z_{\lambda_{n}}(t)\exp\left(
2\frac{z_{\lambda_{n}}^{2}(t)}{z_{0}^{2}}\right)  \mathrm{d}t\\
&  =\frac{z_{0}^{2}}{4}\left(  \exp\left(  2\frac{z_{\lambda_{n}}^{2}%
(T)}{z_{0}^{2}}\right)  -\exp\left(  2\frac{z_{\lambda_{n}}^{2}(0)}{z_{0}^{2}%
}\right)  \right) \\
&  =\left(  1-\lambda_{n}\right)  \int_{0}^{T}p(z_{\lambda_{n}},z_{\lambda
_{n}})(t)z_{\lambda_{n}}(t)\exp\left(  2\frac{z_{\lambda_{n}}^{2}(t)}%
{z_{0}^{2}}\right)  \mathrm{d}t\\
&  -kF_{0}\lambda_{n}\int_{0}^{T}\sin\left(  2kz_{\lambda_{n}}\left(
t\right)  -bt\right)  z_{\lambda_{n}}(t)dt-\frac{4\lambda_{n}F_{0}}{z_{0}^{2}%
}\int_{0}^{T}\left(  \cos\left(  kz_{\lambda_{n}}\left(  t\right)  -\frac
{bt}{2}\right)  \right)  ^{2}z_{\lambda_{n}}^{2}(t)\mathrm{d}t.
\end{align*}

There exist $z_{\lambda_{n}}^{\ast}=z_{\lambda_{n}}(t_{n}^{\ast}%
),\ \widetilde{z}_{\lambda_{n}}=z_{\lambda_{n}}(\widetilde{t}_{n}%
),\ \overline{z}_{\lambda_{n}}=z_{\lambda_{n}}(\overline{t}_{n})$ such that
\begin{align*}
&  (1-\lambda_{n})(z_{\lambda_{n}}^{\ast})^{2}\exp\left(  2\frac{\left(
z_{\lambda_{n}}^{\ast}\right)  ^{2}}{z_{0}^{2}}\right) \\
&  =-kF_{0}\lambda_{n}\sin\left(  2k\widetilde{z}_{\lambda_{n}}-b\widetilde{t}%
_{n}\right)  \widetilde{z}_{\lambda_{n}}-\frac{4\lambda_{n}F_{0}}{Tz_{0}^{2}%
}(\overline{z}_{\lambda_{n}})^{2}\int_{0}^{T}\left(  \cos\left(
kz_{\lambda_{n}}\left(  t\right)  -\frac{bt}{2}\right)  \right)
^{2}\mathrm{d}t,
\end{align*}
so%
\begin{align}
&  \frac{(1-\lambda_{n})}{\lambda_{n}}\left(  \frac{z_{\lambda_{n}}^{\ast}%
}{\overline{z}_{\lambda_{n}}}\right)  ^{2}\exp\left(  2\frac{\left(
z_{\lambda_{n}}^{\ast}\right)  ^{2}}{z_{0}^{2}}\right)
\label{EqContradiction2}\\
&  =-kF_{0}\sin\left(  2k\widetilde{z}_{\lambda_{n}}-b\widetilde{t}%
_{n}\right)  \frac{\widetilde{z}_{\lambda_{n}}}{(\overline{z}_{\lambda_{n}%
})^{2}}-\frac{4F_{0}}{Tz_{0}^{2}}\int_{0}^{T}\left(  \cos\left(
kz_{\lambda_{n}}\left(  t\right)  -\frac{bt}{2}\right)  \right)
^{2}\mathrm{d}t.\nonumber
\end{align}

From \eqref{BVProblemAux} and \eqref{EqIntp} we obtain that%
\begin{align*}
\left\vert z_{\lambda_{n}}(t)-z_{\lambda_{n}}(0)\right\vert  &  \leq\int%
_{0}^{T}\left\vert (1-\lambda_{n})p(z_{\lambda_{n}},z_{\lambda_{n}%
})(s)+\lambda_{n}F_{z}(s,z_{\lambda_{n}}(s))\right\vert \mathrm{d}s\\
&  \leq2\int_{0}^{T}\left\vert F_{z}(s,z_{\lambda_{n}}(s))\right\vert
\mathrm{d}s\\
&  \leq2\int_{0}^{T}\left(  kF_{0}\left\vert f\left(  z_{\lambda_{n}%
}(s)\right)  \right\vert +F_{0}\left\vert f^{\prime}\left(  z_{\lambda_{n}%
}(s)\right)  \right\vert \right)  \mathrm{d}s\\
&  =2F_{0}T(k\left\vert f\left(  \widehat{z}_{\lambda_{n}}\right)  \right\vert
+\left\vert f^{\prime}(\widehat{z}_{\lambda_{n}})\right\vert )\\
&  =2F_{0}T\left(  k\exp\left(  -2\frac{\widehat{z}_{\lambda_{n}}^{2}}%
{z_{0}^{2}}\right)  +4\frac{\widehat{z}_{\lambda_{n}}}{z_{0}^{2}}\exp\left(
-2\frac{\widehat{z}_{\lambda_{n}}^{2}}{z_{0}^{2}}\right)  \right)  ,\text{
}\forall t\in\lbrack0,T],
\end{align*}
for some $\widehat{z}_{\lambda_{n}}=z_{\lambda_{n}}(\widehat{t}_{n})$. This
implies, in particular, that%
\begin{equation}
\sup_{t,s\in\lbrack0,T]}\left\vert z_{\lambda_{n}}(t)-z_{\lambda_{n}%
}(s)\right\vert \rightarrow0,\ \text{as }n\rightarrow\infty. \label{Convergz2}%
\end{equation}

It follows from \eqref{Convergz2} that
\[
\frac{\widetilde{z}_{\lambda_{n}}}{\overline{z}_{\lambda_{n}}}\rightarrow1.
\]
Hence,%
\[
-kF_{0}\sin\left(  2k\widetilde{z}_{\lambda_{n}}-b\widetilde{t}_{n}\right)
\frac{\widetilde{z}_{\lambda_{n}}}{(\overline{z}_{\lambda_{n}})^{2}%
}\rightarrow0,
\]
so \eqref{EqContradiction2} and \eqref{Converg2} (which can be proved in the
same way as in Theorem \ref{ThCotaZetaLambda}) imply that
\[
\frac{(1-\lambda_{n})}{\lambda_{n}}\left(  \frac{z_{\lambda_{n}}^{\ast}%
}{\overline{z}_{\lambda_{n}}}\right)  ^{2}\exp\left(  2\frac{\left(
z_{\lambda_{n}}^{\ast}\right)  ^{2}}{z_{0}^{2}}\right)  \rightarrow
-\frac{2F_{0}}{z_{0}^{2}},
\]
which is not possible as $\frac{(1-\lambda_{n})}{\lambda_{n}}\left(
\frac{z_{\lambda_{n}}^{\ast}}{\overline{z}_{\lambda_{n}}}\right)  ^{2}%
\exp\left(  2\frac{\left(  z_{\lambda_{n}}^{\ast}\right)  ^{2}}{z_{0}^{2}%
}\right)  >0$ for all $n$. Thus, \eqref{EstRoPrima2} is proved.

Finally, \eqref{CotaZLambda2} is deduced exactly in the same way as in Theorem
\ref{ThCotaZetaLambda}.
\end{proof}

\bigskip

Arguing as in Theorem \ref{ThPerPolynomial} we obtain the existence of
periodic solutions for equation \eqref{Equation}.

\begin{theorem}
The boundary-value problem \eqref{Equation}, \eqref{BC} possesses at least one
solution, which can be extended to a periodic solution of equation \eqref{Equation}.
\end{theorem}

\bigskip

Finally, using Theorem \ref{ThBoundPer} we prove that the set of all periodic
solutions of equation \eqref{Equation} is bounded.

\begin{theorem}
There exists a constant $D>0$ such that every solution $z($\textperiodcentered
$)$ of problem \eqref{Equation}, \eqref{BC} satisfies%
\[
\left\Vert z\right\Vert _{C}\leq D.
\]

\end{theorem}

\begin{proof}
It is sufficient to check that the assumptions of Theorem \ref{ThBoundPer} are
satisfied for $f(z)=\exp\left(  -2\frac{z^{2}}{z_{0}^{2}}\right)  $.

Conditions 1--3 are straightforward to verify. Since $f(z)\not =0$ for all
$z$, condition 5 follows from Lemma \ref{FixedCriterion}.

It remains to check the fourth one. For sequences $z_{n},v_{n}$ such that
$z_{n},v_{n}\rightarrow\infty$, $\frac{v_{n}}{z_{n}}\rightarrow1$, as
$n\rightarrow\infty$, and $f^{\prime}\left(  z_{n}\right)  \not =0$, we have%

\begin{align*}
\frac{f^{2}(v_{n})}{f^{\prime}(z_{n})}  &  =\frac{z_{0}^{2}\exp\left(
-\frac{4v_{n}^{2}}{z_{0}^{2}}\right)  }{(-4z_{n})\exp\left(  -2\frac{z_{n}%
^{2}}{z_{0}^{2}}\right)  }=-\frac{z_{0}^{2}}{4z_{n}}\exp\left(  -\frac
{4v_{n}^{2}}{z_{0}^{2}}+2\frac{z_{n}^{2}}{z_{0}^{2}}\right) \\
&  =-\frac{z_{0}^{2}}{4z_{n}}\exp\left(  -\frac{2v_{n}^{2}}{z_{0}^{2}}\left(
2-\frac{z_{n}^{2}}{v_{n}^{2}}\right)  \right)  \rightarrow0,
\end{align*}%
\[
\frac{(f^{\prime}(v_{n}))^{2}}{f^{\prime}(z_{n})}=\frac{16v_{n}^{2}\exp\left(
-\frac{4v_{n}^{2}}{z_{0}^{2}}\right)  }{z_{0}^{2}(-4z_{n})\exp\left(
-\frac{2z_{n}^{2}}{z_{0}^{2}}\right)  }=-\frac{4}{z_{0}^{2}}\frac{v_{n}}%
{z_{n}}v_{n}\exp\left(  -\frac{2v_{n}^{2}}{z_{0}^{2}}\left(  2-\frac{z_{n}%
^{2}}{v_{n}^{2}}\right)  \right)  \rightarrow0.
\]

\end{proof}

\bigskip

\begin{corollary}
Every periodic solution with period $T$ belongs to the ball of radius $D$
centered at $0$ of the space of continuous bounded functions $C_{\mathrm{b}%
}(\mathbb{R}).$
\end{corollary}

\section{Plane waves $f(z)=1$}

If we take $f(z)=1$ in problem \eqref{problem}, the electromagnetic fields
given by equation \eqref{campo} correspond to two dephased counter-propagating
plane waves as we has previously mentioned. This is a limit case when
$z_{0}\rightarrow\infty$ of the previously studied Gaussian and Lorentzian functions.

For plane waves, the problem \eqref{problem} have an analytical solution given
by:
\begin{equation}
z(t)=\frac{b\,t}{2k}-\frac{1}{k}\arctan(H(t)), \label{solu}%
\end{equation}
where
\begin{equation}
H(t)=\frac{2F_{0}k^{2}+\sqrt{b^{2}-4F_{0}^{2}k^{4}}\,\tan\left(  \frac
{\sqrt{b^{2}-4F_{0}^{2}k^{4}}\,\,t}{2}-\arctan(\frac{2F_{0}k^{2}%
+b\,tan(z_{i}k)}{\sqrt{b^{2}-4F_{0}^{2}k^{4}}})\right)  }{b} \label{h}%
\end{equation}
and we have assumed the initial condition $z(0)=z_{i}$.

From this solution it can be deduced that the particles trapped by the
intensity gradient generated by the plane waves are translated indefinitely
with nearly the constant velocity phase of the conveyor $v_{c}=\frac{b}{2k}$
\cite{Ruffner2012}.

Depending on the problem parameters two situations can be observed. On the one
hand, if $b\leq2F_{0}k^{2}$ the particles trajectories correspond to a uniform
rectilinear motion. On the other hand, if $b>2F_{0}k^{2}$, the particles
describe a bounded oscillating path around the straight line $z=z_{i}%
+\frac{b\,t}{2k}-\frac{1}{2k}\sqrt{b^{2}-4F_{0}^{2}k^{4}}t$. In both cases, as
it has been previously mentioned, the particles are indefinitely translated,
but plane waves have not enough intensity to move particles and other beams
($f(z)\neq1$) must be used.

We observe that the map $\arctan$ is multivalued and when $b>2F_{0}k^{2}$ we
need to choose the branch of the function in such a way that the identity
$\arctan(\tan(x))=x$ is satisfied.

Moreover, this limit case can be used too to show that the solutions to
problem \eqref{problem} depend strongly on the parameters of the problem. For
example, an approximated periodic solution can be found if $F_{0}\approx0$ in
the case of plane waves ($f(z)=1$). Using a Taylor series expansion of order
one of function $H(t)$ (given by equation \eqref{h}) around $F_{0}<<1$, and
introducing it in solution \eqref{solu}, we obtain that%
\begin{equation}
z(t)\approx z_{i}+\frac{F_{0}\,k\,(\cos(2\,z_{i}\,k)-\cos(2\,z_{i}%
\,k-b\,t))}{b},
\end{equation}
which shows that for low values of the parameter $F_{0}$ the particles
describes approximately a periodic movement around its initial position.

\section{Numerical simulations}

In the next section we show the results obtained for the parameters
$F_{0}=0.8\,pm/s$, $b=100\,Hz$, $z_{0}=0.37\,\lambda$, $k=2.66\,\pi/\lambda$
and $\lambda=580\,nm$.

\begin{figure}[tbh]
\centerline{\scalebox{0.5}{\includegraphics[angle=0]{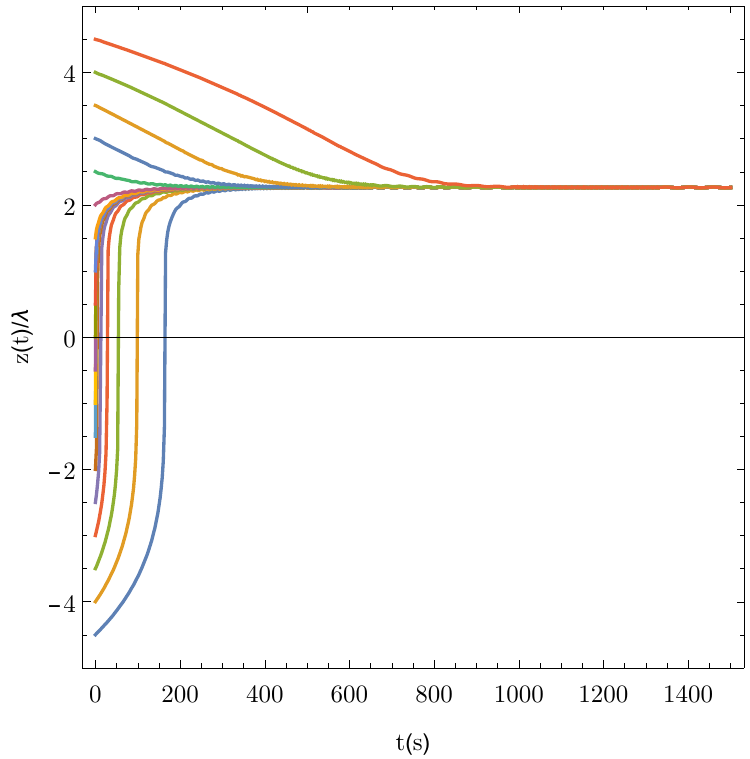}}}\caption{{\protect\footnotesize {Particle's
trajectories for different initial conditions for $f(z)=\frac{1}%
{1+(z/z_{0})^{2}}$.}}}%
\label{fig1}%
\end{figure}

\begin{figure}[tbh]
\centerline{\scalebox{0.5}{\includegraphics[angle=0]{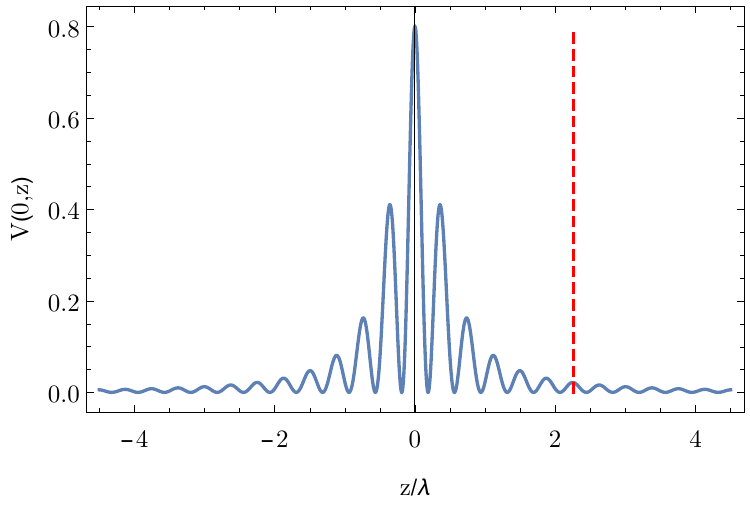}}}\caption{{\protect\footnotesize {Potential
V(0,z) for $f(z)=\frac{1}{1+(z/z_{0})^{2}}$. The dotted line shows the z
convergence position.}}}%
\label{pot1}%
\end{figure}

\begin{figure}[tbh]
\centerline{\scalebox{0.5}{\includegraphics[angle=0]{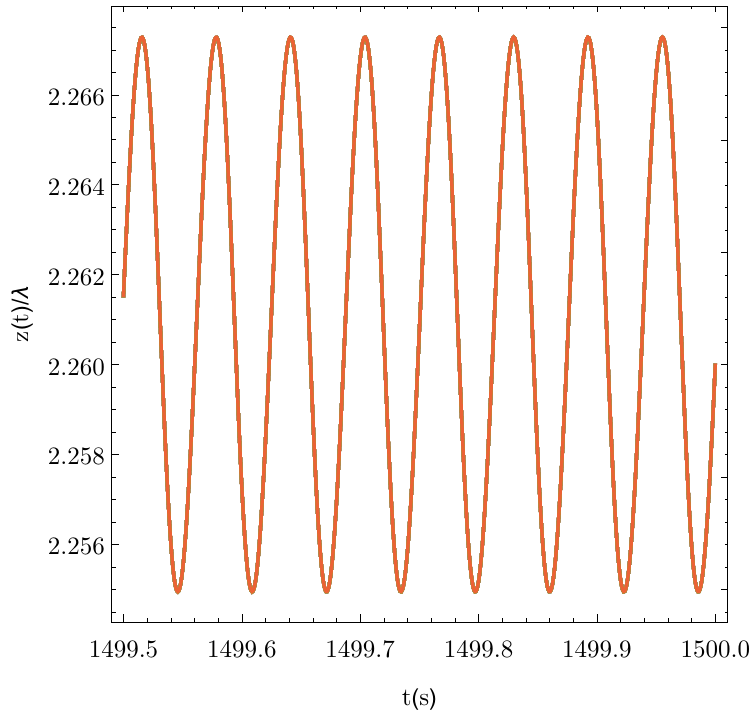}}}\caption{{\protect\footnotesize {Particle's
trajectories behavior for different initial conditions at convergence region
for $f(z)=\frac{1}{1+(z/z_{0})^{2}}$.}}}%
\label{fig2}%
\end{figure}

For the case of counter-propagating Gaussian beams ($f(z)=\frac{1}%
{1+(z/z_{0})^{2}}$) theoretically analyzed in previous sections the
trajectories for the given parameters are shown in Fig. \ref{fig1}. As it can
be seen, all particles initially located in the interval $[-4.5\lambda
,4.5\lambda]$ (where the potential function $V(t,z)$ shows significant values,
see Fig. \ref{pot1}) converge toward the same $z$ region, showing all of them,
as it can be observed in Fig. \ref{fig2}, a periodic behavior with the same
frequency and amplitude.

\newpage

\begin{figure}[tbh]
\centerline{\scalebox{0.5}{\includegraphics[angle=0]{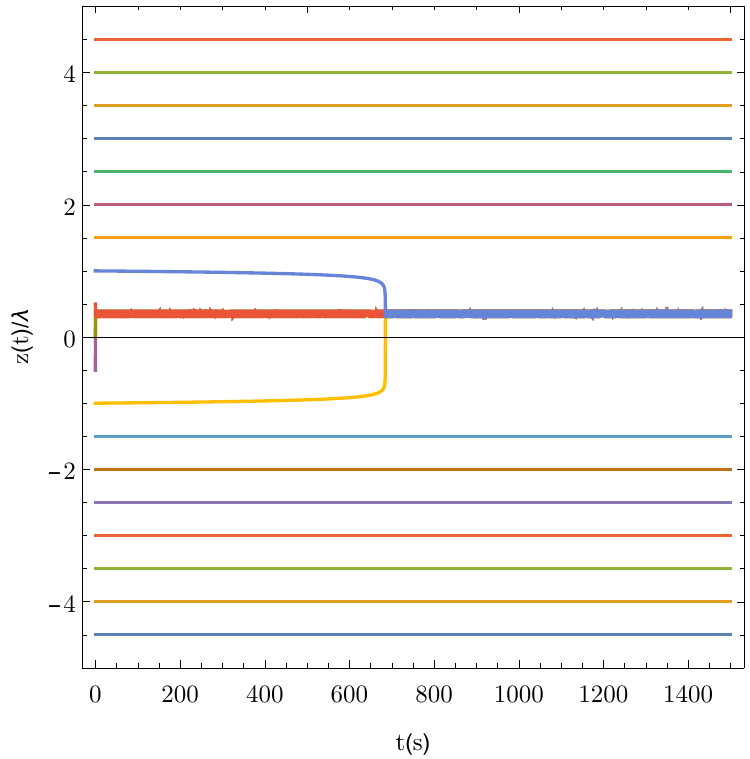}}}\caption{{\protect\footnotesize {Particle's
trajectories behavior for different initial conditions at convergence region
for $f(z)=exp(-2\frac{z^{2}}{z_{0}^{2}})$.}}}%
\label{fig3}%
\end{figure}

\begin{figure}[tbh]
\centerline{\scalebox{0.5}{\includegraphics[angle=0]{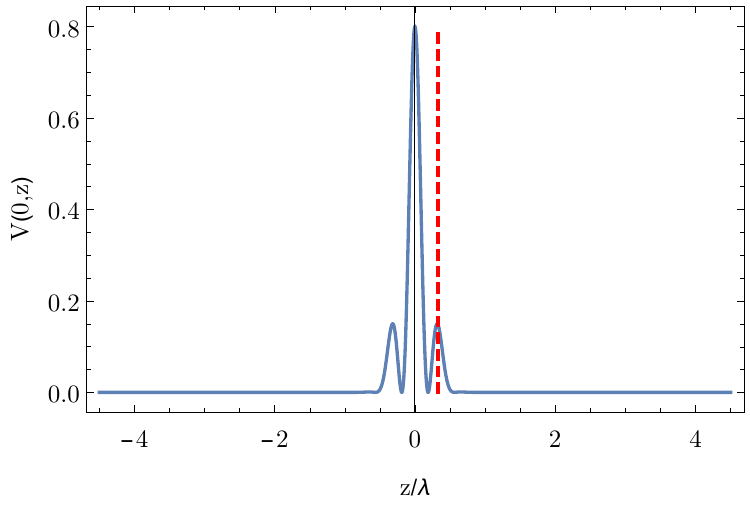}}}\caption{{\protect\footnotesize {V(0,z)
potential for $f(z)=\exp(-2\frac{z^{2}}{z_{0}^{2}})$. The dotted line shows
the z convergence position.}}}%
\label{pot2}%
\end{figure}

\begin{figure}[tbh]
\centerline{\scalebox{0.5}{\includegraphics[angle=0]{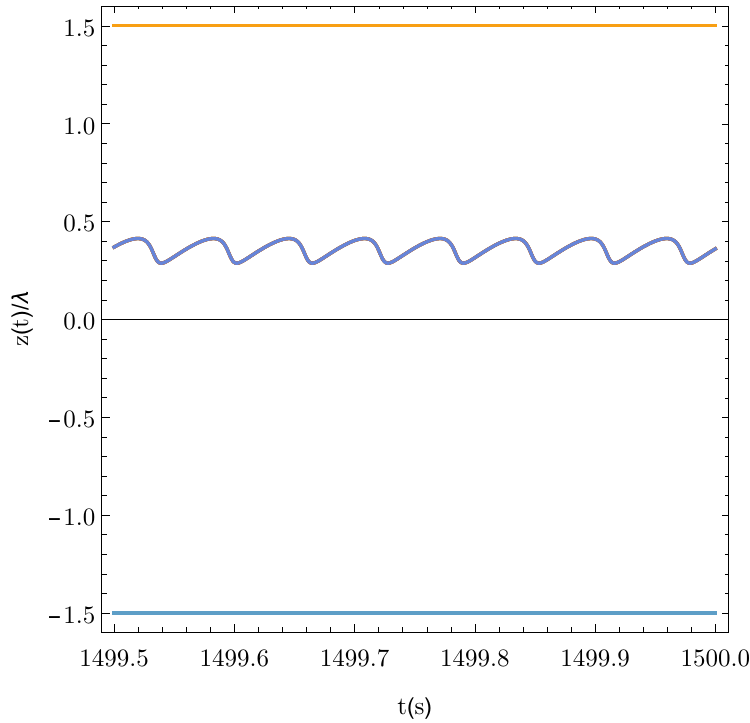}}}\caption{{\protect\footnotesize {Particle's
trajectories behavior for different initial conditions at convergence region
for $f(z)=\exp(-2\frac{z^{2}}{z_{0}^{2}})$.}}}%
\label{fig4}%
\end{figure}

Fig. \ref{fig3} shows the particle's trajectories of problem \eqref{problem}
when $f(z)=\exp(-2\frac{z^{2}}{z_{0}^{2}})$. As it can be observed in this
case, only the particles initially located at interval $[-\lambda,\lambda]$
converge toward the same region, showing again all of them a periodic behavior
with the same frequency and amplitude (see Fig. \ref{fig4}). The particles
outside of the interval $[-\lambda,\lambda]$ remain in this case at the
initial position because the potential energy is null in this region (see Fig.
\ref{pot2}).

In the limit case $z_{0}\rightarrow\infty$ of plane waves ($f(z)=1$) and for
the same value parameters, the conveyor speed is $v_{c}=\frac{b}%
{2k}=5.98\lambda\,nm/s$, so, according to solution \eqref{solu} and taking
into account that $b<2F_{0}k^{2}$, the particles will describe a rectilinear
motion $\frac{z(t)}{\lambda}\approx z_{i}+5.98\,t$ which implies that
$\frac{z(1500)}{\lambda}\approx9000$ and there is not an oscillating behavior.

%\bigskip
\pagebreak

\textbf{Acknowledgements}

The second author was partially supported by Spanish Ministry of Economy and
Competitiveness and FEDER, projects MTM2015-63723-P and MTM2016-74921-P, and
by Junta de Andaluc\'{\i}a (Spain), project P12-FQM-1492.

We would like to thank the anonymous referees and also the associate editor
for the careful reading of the manuscript and their helpful suggestions, which
allowed us to improve the paper.

%\pagebreak

%\bibliography{/home/luis/disco/BIBLIOGRAFIA/references}

%\bibliographystyle{/home/luis/disco/BIBLIOGRAFIA/osajnl.bst}

\end{document}